\documentstyle[12pt]{article}

\newcommand{\be}{\begin{equation}}
\newcommand{\ee}{\end{equation}}
\newcommand{\ba}{\begin{eqnarray}}
\newcommand{\ea}{\end{eqnarray}}
\newcommand{\ban}{\begin{eqnarray*}}
\newcommand{\ean}{\end{eqnarray*}}
\newcommand{\n}{\nonumber \\}

\newcommand{\eq}[1]{(\ref{#1})}
\newcommand{\sfrac}[2]{{\textstyle \frac{#1}{#2}}}
\newcommand{\ignore}[1]{}

\newtheorem{prop}{Proposition}
\newcommand{\qed}{\hbox{\rule[-1pt]{4pt}{9pt}}}

\font\csc=cmcsc10 scaled\magstep1
\font\tennn=msbm10
\font\twelvenn=msbm10 scaled\magstep1

\newcommand{\Sbm}[1]{\leavevmode\raise-.15ex\hbox{\twelvenn #1}}
\newcommand{\sbm}[1]{\leavevmode\raise-.15ex\hbox{\tennn #1}}

\newcommand{\bZ}{\Sbm{Z}}
\newcommand{\bz}{\sbm{Z}}

\begin{document}

\renewcommand{\thefootnote}{\fnsymbol{footnote}}
\font\csc=cmcsc10 scaled\magstep1

{\baselineskip=14pt
 \rightline{
 \vbox{
       \hbox{EFI-97-07}
       \hbox{DPSU-97-1}
       \hbox{January 1997}
}}}

\vskip 11mm
\begin{center}

{\large\bf $q$-Difference Realization of $U_q(sl(M|N))$ \\
and \\\vskip5pt
Its Application to Free Boson Realization of $U_q(\widehat{sl}(2|1))$
}

\vspace{15mm}

{\csc Hidetoshi AWATA}\footnote[1]{JSPS fellow}$^{1}$,
{\csc Satoru ODAKE}$\,{}^2$ and
{\csc Jun'ichi SHIRAISHI}$\,{}^3$
\\ 
%
{\baselineskip=15pt
\it\vskip.35in 
\setcounter{footnote}{0}\renewcommand{\thefootnote}{\arabic{footnote}}
\footnote{e-mail address : awata@rainbow.uchicago.edu}
James Frank Institute and Enrico Fermi Institute,
University of Chicago,\\
5640 S. Ellis Ave., Chicago, IL 60637, U.S.A.
\vskip.1in 
\footnote{e-mail address : odake@azusa.shinshu-u.ac.jp}
Department of Physics, Faculty of Science \\
Shinshu University, Matsumoto 390, Japan\\
\vskip.1in 
\footnote{e-mail address : shiraish@momo.issp.u-tokyo.ac.jp}
Institute for Solid State Physics, \\
University of Tokyo, Tokyo 106, Japan \\
}
\end{center}

\vspace{7mm}
\begin{abstract}
We present a $q$-difference realization of the quantum superalgebra 
$U_q(sl(M|N))$, which includes Grassmann even and odd coordinates
and their derivatives.
Based on this result we obtain a free boson realization of the quantum
affine superalgebra $U_q(\widehat{sl}(2|1))$ of an arbitrary level $k$.
\end{abstract}

\vspace{30mm}
q-alg/9701032

\vfill\eject
\setcounter{footnote}{0}\renewcommand{\thefootnote}{\arabic{footnote}}

\section{Introduction}

Free field realization is the most useful tool for studying representation
theories of infinite dimensional algebras and their application to
physics such as calculation of correlation functions.
The Virasoro and affine Lie algebras are investigated extensively by this
method\cite{rBMP}.
Quantum deformation ($q$-deformation) of these algebras are also
studied \cite{rSKAOAKOSFFr,rLP,rDFJMNJMMN,rAOS}.
The aim of this article and its continuation\cite{rAOSU} 
is to generalize these to superalgebra cases,
which will help us to study physical systems with superalgebra symmetries
such as supersymmetric $t$-$J$ model, to construct a supersymmetric
quantum deformed Virasoro algebra, etc.
A free boson realization of $U_q(\widehat{sl}(M|N))$ has been obtained
for a level one\cite{rKSU} but not yet for an arbitrary level $k$.
To obtain a free boson realization of $U_q(\widehat{sl}(M|N))$ with
an arbitrary level (Wakimoto module), we first consider a $q$-difference
realization of $U_q(sl(M|N))$, because for non-superalgebra case 
the $q$-difference realization of $U_q(sl_N)$ 
on the flag manifold\cite{rANO} helped us
to construct the free boson realization of $U_q(\widehat{sl}_N)$ \cite{rAOS}.
Next we try to construct a free boson realization of $U_q(\widehat{sl}(M|N))$.
At present we have obtained it for $M+N\leq 3$. 
In this article the result for $(M,N)=(2,1)$ is given for simplicity.

After recalling the definitions of $U_q(sl(M|N))$ and 
$U_q(\widehat{sl}(M|N))$ in section 2, we present a $q$-difference
realization of $U_q(sl(M|N))$ in section 3.
In section 4 we give a free boson realization of $U_q(\widehat{sl}(2|1))$ 
with an arbitrary level $k$.

While preparing our draft, we received ref.\cite{rK}, where
Kimura obtained a $q$-difference realization of $U_q(sl(M|N))$.
We will comment on it at the end of section 3.

\section{$U_q(sl(M|N))$ and $U_q(\widehat{sl}(M|N))$}

First let us recall the definition of $U_q(sl(M|N))$ 
($M,N,M+N-2\geq 0$)\cite{rY}.
We take the Cartan matrix of $sl(M|N)$ as
$a_{ij}=(\nu_i+\nu_{i+1})\delta_{i,j}
-\nu_i\delta_{i,j+1}-\nu_{i+1}\delta_{i+1,j}$ ($1\leq i,j\leq M+N-1$),
where $\nu_i$ 
is $1$ for $1\leq i\leq M$, $-1$ for $M+1\leq i\leq M+N$. 
The quantum Lie superalgebra $U_q(sl(M|N))$ is defined by the Chevalley
generators $t_i=q^{h_i}$ (invertible), $e^+_i=e_i$, $e^-_i=f_i$ 
($i=1,\cdots,M+N-1$) with the relations,
\ba
  &&
  t_it_j=t_jt_i \quad([h_i,h_j]=0), \\
  &&
  t_ie^{\pm}_jt_i^{-1}=q^{\pm a_{ij}}e^{\pm}_j \quad
  ([h_i,e^{\pm}_j]=\pm a_{ij}e^{\pm}_j),
  \label{he} \\
  &&
  [e_i,f_j]=\delta_{i,j}[h_i],
  \label{ef} \\
  &&
  [e^{\pm}_i,[e^{\pm}_i,e^{\pm}_j]_{q^{-1}}]_q=0 \quad
  \mbox{ for } |a_{ij}|=1, i\neq M,\\
  &&
  [e^{\pm}_M,[e^{\pm}_{M+1},[e^{\pm}_M,e^{\pm}_{M-1}]_{q^{-1}}]_q]=0,
\ea
where $[x]=(q^x-q^{-x})/(q-q^{-1})$,
$[A,B]_{\xi}=AB-(-1)^{|A|\cdot|B|}\xi BA$, and $[A,B]=[A,B]_1$.
Here $|\cdot|$ stands for $\bZ_2$-grading (Grassmann parity)
($|e^{\pm}_M|=1$ and 0 for otherwise).

Next recall the definition of the quantum affine superalgebra 
$U_q(\widehat{sl}(M|N))$ in terms of Drinfeld generators, 
$E^{+,i}_n=E^i_n$, $E^{-,i}_n=F^i_n$ ($n\in\bZ$), 
$H^i_n$ ($n\in\bZ_{\neq 0}$), invertible $K_i$ ($1\leq i\leq M+N-1$) 
and invertible $\gamma$ \cite{rY}.
Their $\bZ_2$ gradings are $1$ for $E^{\pm,M}_n$ and zero otherwise.
The relations are
\ba
  &&
  \gamma \mbox{ : central element},\\
  &&
  [K_i,H^j_n]=0,\quad
  K_iE^{\pm,j}_nK_i^{-1}=q^{\pm a_{ij}}E^{\pm,j}_n,\\
  &&
  [H^i_n,H^j_m]=
  \frac{1}{n}[a_{ij}n]\frac{\gamma^n-\gamma^{-n}}{q-q^{-1}}
   \delta_{n+m,0},\\
  &&
  [H^i_n,E^{\pm,j}_m]=
  \pm\frac{1}{n}[a_{ij}n]\gamma^{\mp\frac12|n|}E^{\pm,j}_{n+m},\\
  &&
  [E^{+,i}_n,E^{-,j}_m]=
  \frac{\delta^{ij}}{q-q^{-1}}\Bigl(
  \gamma^{\frac12(n-m)}\psi^i_{+,n+m}-\gamma^{-\frac12(n-m)}\psi^i_{-,n+m}
  \Bigr),\\
  &&
  [E^{\pm,i}_{n+1},E^{\pm,j}_m]_{q^{\pm a_{ij}}}+
  [E^{\pm,j}_{m+1},E^{\pm,i}_n]_{q^{\pm a_{ij}}}=0,\\
  &&
  [E^{\pm,i}_n,E^{\pm,j}_m]=0\quad \mbox{ for }a_{ij}=0,\\
  &&
   [E^{\pm,i}_{n_1},[E^{\pm,i}_{n_2},E^{\pm,j}_m]_{q^{-1}}]_q
  +(n_1\leftrightarrow n_2)
  =0 \quad \mbox{ for } |a_{ij}|=1,i\neq M,\\
  &&
   [E^{\pm,M}_{n_1},[E^{\pm,M+1}_m,[E^{\pm,M}_{n_2},
                                    E^{\pm,M-1}_{\ell}]_{q^{-1}}]_q]
  +(n_1\leftrightarrow n_2)
  =0,
\ea
where $\psi^i_{\pm,n}$ is defined by
\be
  \sum_{n\in\bz}\psi^i_{\pm,n}z^{-n}=
  K_i^{\pm 1}\exp\Bigl(\pm(q-q^{-1})\sum_{\pm n>0}H^i_nz^{-n}\Bigr).
\ee
Let $H^i_0$ be defined by
\be
  K_i=\exp\Bigl((q-q^{-1})\sfrac{1}{2}H^i_0\Bigr),
\ee
then above relations hold for $H^i_n$ ($n\in\bZ$)
(in the case of $n=0$, $\frac{1}{n}\ast$ should be understood as 
$\lim_{n\rightarrow 0}\frac{1}{n}\ast$).
Fields $H^i(z)$, $E^{\pm,i}(z)$ and $\psi_{\pm}^i(z)$ are defined by
\be
  H^i(z)=\sum_{n\in\bz}H^i_nz^{-n-1},\quad
  E^{\pm,i}(z)=\sum_{n\in\bz}E^{\pm,i}_nz^{-n-1},\quad
  \psi_{\pm}^i(z)=\sum_{n\in\bz}\psi_{\pm,n}^iz^{-n}.
\ee

\section{$q$-Difference realization of $U_q(sl(M|N))$}

In ref.\cite{rANO} $U_q(sl_N)$ was realized by $q$-difference operators
in flag coordinates. We extend it to the super case.

Let $x_{i,j}$ ($1\leq i<j\leq M+N$) be variables with Grassmann parity 
$\nu_i\nu_j$ and 
set $\vartheta_{i,j}=x_{i,j}\frac{\partial}{\partial x_{i,j}}$.
For $1\leq i,j\leq M+N-1$, we define $h_i$, $e_{i,i}$, 
$e_{i,i'}$ ($1\leq i'\leq i-1$), $f_{j,j'}^1$ ($1\leq j'\leq j-1$), 
$f_{j,j}^2$ and $f_{j,j'}^3$ ($j+2\leq j'\leq M+N$) as follows:
\ba
  h_i 
  &\!\!=\!\!&
  -\sum_{j=1}^{i-1}(\nu_{i+1}\vartheta_{j,i+1}-\nu_i\vartheta_{j,i})
  +\lambda_i-(\nu_i+\nu_{i+1})\vartheta_{i,i+1}
  -\sum_{j=i+2}^{M+N}(\nu_i\vartheta_{i,j}-\nu_{i+1}\vartheta_{i+1,j}),\\
  e_{i,i}
  &\!\!=\!\!&
  \frac{1}{x_{i,i+1}}[\vartheta_{i,i+1}]
  q^{-\sum_{\ell=1}^{i-1}
  (\nu_{i+1}\vartheta_{\ell,i+1}-\nu_i\vartheta_{\ell,i})}, \\
  e_{i,i'}
  &\!\!=\!\!&
  x_{i',i}\frac{1}{x_{i',i+1}}[\vartheta_{i',i+1}]
  q^{-\sum_{\ell=1}^{i'-1}
  (\nu_{i+1}\vartheta_{\ell,i+1}-\nu_i\vartheta_{\ell,i})}, \\
  f_{j,j'}^1
  &\!\!=\!\!&
  x_{j',j+1}\frac{1}{x_{j',j}}[\vartheta_{j',j}] \n
  &&\times
  q^{\sum_{m=j'+1}^{j-1}
      (\nu_{j+1}\vartheta_{m,j+1}-\nu_j\vartheta_{m,j})
     -\lambda_j+(\nu_j+\nu_{j+1})\vartheta_{j,j+1}
     +\sum_{m=j+2}^{M+N}
      (\nu_j\vartheta_{j,m}-\nu_{j+1}\vartheta_{j+1,m})}, \\
  f_{j,j}^2
  &\!\!=\!\!&
  x_{j,j+1}\Bigl[\lambda_j-\nu_j\vartheta_{j,j+1}
   -\sum_{m=j+2}^{M+N}
    (\nu_j\vartheta_{j,m}-\nu_{j+1}\vartheta_{j+1,m})\Bigr], \\
  f_{j,j'}^3
  &\!\!=\!\!&
  x_{j,j'}\frac{1}{x_{j+1,j'}}[\vartheta_{j+1,j'}]
  q^{\lambda_j
     -\sum_{m=j'}^{M+N}
      (\nu_j\vartheta_{j,m}-\nu_{j+1}\vartheta_{j+1,m})},
\ea
where, for Grassmann odd variable $x$, 
the expression $\frac{1}{x}$ stands for the derivative 
$\frac{1}{x}=\frac{\partial}{\partial x}$.

\begin{prop}
The following {\rm (\romannumeral1)} and {\rm (\romannumeral2)} realize
$U_q(sl(M|N))$:
\ba
  &{\rm (\romannumeral1)}&
  h_i, \quad
  e_i=e_{i,i}+\nu_i\sum_{i'=1}^{i-1}e_{i,i'},\quad
  f_j=\sum_{j'=1}^{j-1}f_{j,j'}^1+f_{j,j}^2-\sum_{j'=j+2}^{M+N}f_{j,j'}^3, \\
  &{\rm (\romannumeral2)}&
  h_i, \quad
  e_i=e_{i,i}+\sum_{i'=1}^{i-1}e_{i,i'},\quad
  f_j=\nu_j\sum_{j'=1}^{j-1}f_{j,j'}^1+f_{j,j}^2
      -\nu_{j+1}\sum_{j'=j+2}^{M+N}f_{j,j'}^3.
\ea
\end{prop}
{\it Proof:}
Direct calculation shows this proposition.
For reader's convenience we present some intermediate steps;
\ba
  \!\!\!\!\!\!\!\!\!\!&&
  \lbrack e_{i,i},f_{j,j'}^1 \rbrack
  =\lbrack e_{i,i},f_{j,j'}^3 \rbrack
  =\lbrack e_{i,i'},f_{j,j}^2 \rbrack=0, \\
  \!\!\!\!\!\!\!\!\!\!&&
  \lbrack e_{i,i},f_{j,j}^2 \rbrack
  =
  \delta_{i,j}\Bigl[\lambda_i-(\nu_i+\nu_{i+1})\vartheta_{i,i+1}
  -\sum_{m=i+2}^{M+N}(\nu_i\vartheta_{i,m}-\nu_{i+1}\vartheta_{i+1,m})
  \Bigr]q^{-\sum_{\ell=1}^{i-1}
   (\nu_{i+1}\vartheta_{\ell,i+1}-\nu_i\vartheta_{\ell,i})} \n
  \!\!\!\!\!\!\!\!\!\!&&\qquad\quad\quad
  +\delta_{i,j+1}\nu_ix_{i-1,i}\frac{1}{x_{i,i+1}}[\vartheta_{i,i+1}] \\
  \!\!\!\!\!\!\!\!\!\!&&\qquad\qquad\qquad\times
  q^{-\sum_{\ell=1}^{i-1}
     (\nu_{i+1}\vartheta_{\ell,i+1}-\nu_i\vartheta_{\ell,i})
     +\lambda_{i-1}-\nu_{i-1}\vartheta_{i-1,i}
     -\sum_{m=i+1}^{M+N}
     (\nu_{i-1}\vartheta_{i-1,m}-\nu_i\vartheta_{i,m})}, \n
  \!\!\!\!\!\!\!\!\!\!&&
  \lbrack e_{i,i'},f_{j,j'}^1 \rbrack
  =
  \delta_{i,j}\delta_{i',j'}
  \nu_i[\nu_i\vartheta_{i',i}-\nu_{i+1}\vartheta_{i',i+1}]
  q^{-\sum_{\ell=1}^{i'-1}
     (\nu_{i+1}\vartheta_{\ell,i+1}-\nu_i\vartheta_{\ell,i})
     +\sum_{\ell=i'+1}^{i-1}
     (\nu_{i+1}\vartheta_{\ell,i+1}-\nu_i\vartheta_{\ell,i})} \n
  \!\!\!\!\!\!\!\!\!\!&&\qquad\qquad\qquad\times
  q^{-\lambda_i+(\nu_i+\nu_{i+1})\vartheta_{i,i+1}
     +\sum_{\ell=i+2}^{M+N}
     (\nu_i\vartheta_{i,\ell}-\nu_{i+1}\vartheta_{i+1,\ell})}, \\
  \!\!\!\!\!\!\!\!\!\!&&
  \lbrack e_{i,i'},f_{j,j'}^3 \rbrack
  =
  x_{ji}\frac{1}{x_{j+1,i+1}}[\vartheta_{j+1,i+1}]
  \Bigl(\delta_{i',j}\delta_{j',i+1}-\delta_{i',j+1}\delta_{j',i}
   q^{(\nu_i-\nu_j)\vartheta_{j,i}}\Bigr) \n
  \!\!\!\!\!\!\!\!\!\!&&\qquad\qquad\qquad\times
  q^{-\sum_{\ell=1}^{j-1}
     (\nu_{i+1}\vartheta_{\ell,i+1}-\nu_i\vartheta_{\ell,i})
     +\lambda_j-\nu_{i+1}\vartheta_{j,i+1}
     -\sum_{m=i+1}^{M+N}
     (\nu_j\vartheta_{j,m}-\nu_{j+1}\vartheta_{j+1,m})}.
\ea
Using these results and the formula
\be
  [a]q^{\sum_{i=1}^nb_i}
  +\sum_{i=1}^n[b_i]q^{-a+\sum_{j=1}^{i-1}b_j-\sum_{j=i+1}^nb_j}
  =[a+\sum_{i=1}^nb_i],
\ee
we have
\be
  [e_{i,i},f_{j,j}^2]
  +\epsilon\sum_{i'=1}^{i-1}\sum_{j'=1}^{j-1}[e_{i,i'},f_{j,j'}^1]
  -\epsilon'\sum_{i'=1}^{i-1}\sum_{j'=j+2}^{M+N}[e_{i,i'},f_{j,j'}^3]
  =\delta_{i,j}[h_i],
  \label{efh}
\ee
where $\epsilon=\nu_i,\nu_j$ and $\epsilon'=\nu_i,\nu_{j+1}$ 
(four choices). \hfill \qed

\noindent
{\bf Remark 1~}
(\romannumeral2) is obtained from (\romannumeral1) by replacement 
$x_{i,j}\rightarrow\nu_jx_{i,j}$ and 
$e^{\pm}_i\rightarrow\nu_{i+1}e^{\pm}_i$.

\noindent
{\bf Remark 2~}
Recently Kimura obtained a similar result\cite{rK}.
He uses the following $f_{j,j}^{\prime\,2}$ instead of $f_{j,j}^2$,
\be
  f_{j,j}^{\prime\,2}
  =
  x_{j,j+1}\Bigl[\lambda_j-\sfrac{1}{2}(\nu_j+\nu_{j+1})\vartheta_{j,j+1}
   -\sum_{m=j+2}^{M+N}
    (\nu_j\vartheta_{j,m}-\nu_{j+1}\vartheta_{j+1,m})\Bigr].
\ee
It can be checked that $f_{j,j}^{\prime\,2}=f_{j,j}^2$.
Therefore Proposition 3 in ref.\cite{rK} agrees with our (\romannumeral2)
(Note that $x_{ij}^{\rm Kimura}=x_{i,j+1}$ and 
we have replaced $q$ by $q^{-1}$ such that the realization for $N=0$ case, 
$U_q(sl(M|0))=U_q(sl_M)$, reduces to the result in \cite{rANO}).

\section{Free boson realization of $U_q(\widehat{sl}(2|1))$}

In ref.\cite{rAOS}, 
on the bases of a $q$-difference realization of $U_q(sl_N)$,
a free boson realization of $U_q(\widehat{sl}_N)$ was obtained.
Here we try to generalize it to the $U_q(\widehat{sl}(M|N))$.

Let us introduce oscillators and coordinates 
$a^i_n$, $Q_a^i$ ($n\in\bZ$, $1\leq i\leq M+N-1$), 
$b^{ij}_n$, $Q_b^{ij}$ and $c^{ij}_n$, $Q_c^{ij}$ ($n\in\bZ$, 
$1\leq i<j\leq M+N$), which are all Grassmann even operators and
satisfy the following relations,
\ba
  \lbrack a^i_n,a^j_m \rbrack
  &\!\!=\!\!&
  \frac{1}{n}[(k+g)n][a_{ij}n]\delta_{n+m,0},\quad\,
  [a^i_n,Q_a^j]=(k+g)a_{ij}\delta_{n0},\\
  \lbrack b^{ij}_n,b^{i'j'}_m \rbrack
  &\!\!=\!\!&
  -\nu_i\nu_j\frac{1}{n}[n]^2\delta^{ii'}\delta^{jj'}\delta_{n+m,0},\quad
  [b^{ij}_n,Q_b^{i'j'}]=-\nu_i\nu_j\delta^{ii'}\delta^{jj'}\delta_{n0},\\
  \lbrack c^{ij}_n,c^{i'j'}_m \rbrack
  &\!\!=\!\!&
  \frac{1}{n}[n]^2\delta^{ii'}\delta^{jj'}\delta_{n+m,0},\qquad\quad\;\:
  [c^{ij}_n,Q_c^{i'j'}]=\delta^{ii'}\delta^{jj'}\delta_{n0},
\ea
and other commutators vanish. 
Here $g=M-N$ and $k$ is a complex parameter.
For a pair of oscillator and coordinate $(a_n,Q)$, we define boson fields 
$a(z)$ and $a_{\pm}(z)$ as follows:
\ba
  a(z)
  &\!\!=\!\!&
  -\sum_{n\neq 0}\frac{a_n}{[n]}z^{-n}+Q+a_0\log z,\\
  a_{\pm}(z)
  &\!\!=\!\!&
  \pm\Bigl((q-q^{-1})\sum_{\pm n>0}a_nz^{-n}+a_0\log q\Bigr).
\ea
Normal ordering $:~:$ is defined as follows;
move $a_n$ ($n\geq 0$) to right, $a_n$ ($n<0$) and $Q$ to left,
and arrange the zero modes $e^{\pm Q_b^{ij}}$ with $\nu_i\nu_j=-1$
in  lexicographic ordering for $(ij)$.
It is well known in CFT that when we bosonize many fermions we have to
introduce cocycle factors such that different fermions anticommute.
In our case cocycle factors are needed for $:\!e^{\pm b^{ij}(z)}\!:$ with 
$\nu_i\nu_j=-1$. 
But we suppress them by modifying the commutation relation for 
$Q_b^{ij}$ with $\nu_i\nu_j=-1$;
set $[Q_b^{ij},Q_b^{i'j'}]=i\pi$, then we have a minus sign by
interchanging operators $e^{\pm Q_b^{ij}}$, 
$e^{\epsilon Q_b^{ij}}e^{\epsilon' Q_b^{i'j'}}=
-e^{\epsilon' Q_b^{i'j'}}e^{\epsilon Q_b^{ij}}$ 
($\nu_i\nu_j=\nu_{i'}\nu_{j'}=-1$, $(i,j)\neq(i',j')$, 
$\epsilon,\epsilon'=\pm 1$).
For example, for $(M,N)=(2,1)$, 
$:e^{Q_b^{23}}e^{Q_b^{13}}:~=~:e^{Q_b^{13}}e^{Q_b^{23}}:~=
e^{Q_b^{23}}e^{Q_b^{13}}=-e^{Q_b^{13}}e^{Q_b^{23}}$.

In the following we restrict ourselves to $(M,N)=(2,1)$ case.
Let us define fields $H^i(z)$, $E^i(z)$ and $F^i(z)$ 
($i=1,2$) as follows:
\ba
  H^1(z)
  &\!\!=\!\!&
  \frac{1}{(q-q^{-1})z}\Bigl(
  a^1_+(q^{\frac12}z)
  +b^{12}_+(q^{\frac{k}{2}}z)+b^{12}_+(q^{\frac{k}{2}+2}z)
  +b^{13}_+(q^{\frac{k}{2}+2}z)-b^{23}_+(q^{\frac{k}{2}+1}z) \n
  &&\quad
  -a^1_-(q^{-\frac12}z)
  -b^{12}_-(q^{-\frac{k}{2}}z)-b^{12}_-(q^{-\frac{k}{2}-2}z)
  -b^{13}_-(q^{-\frac{k}{2}-2}z)+b^{23}_-(q^{-\frac{k}{2}-1}z)\Bigr),
  \label{H1}\\
  H^2(z)
  &\!\!=\!\!&
  \frac{1}{(q-q^{-1})z}\Bigl(
  a^2_+(q^{\frac12}z)
  -b^{12}_+(q^{\frac{k}{2}+1}z)-b^{13}_+(q^{\frac{k}{2}+1}z) \n
  &&\quad
  -a^2_-(q^{-\frac12}z)
  +b^{12}_-(q^{-\frac{k}{2}-1}z)+b^{13}_-(q^{-\frac{k}{2}-1}z)\Bigr),
  \label{H2}\\
  E^1(z)
  &\!\!=\!\!&
  \frac{-1}{(q-q^{-1})z}:\Bigl(
  e_{11}e^{b^{12}_+(z)-(b+c)^{12}(qz)}
  -e_{12}e^{b^{12}_-(z)-(b+c)^{12}(q^{-1}z)}\Bigl):,
  \label{E1}\\
  E^2(z)
  &\!\!=\!\!&
  :\Bigl(
  e_{21}e^{-b^{12}_+(qz)-b^{13}_+(qz)+b^{23}(qz)}
  +e_{22}e^{(b+c)^{12}(z)+b^{13}(z)}\Bigr):,
  \label{E2}\\
  F^1(z)
  &\!\!=\!\!&
  \frac{1}{(q-q^{-1})z}:\Bigl(
  f_{11}e^{a^1_+(q^{\frac{k+1}{2}}z)+b^{12}_+(q^{k+2}z)
  +b^{13}_+(q^{k+2}z)-b^{23}_+(q^{k+1}z)+(b+c)^{12}(q^{k+1}z)} \n
  &&\qquad
  -f_{12}e^{a^1_-(q^{-\frac{k+1}{2}}z)+b^{12}_-(q^{-k-2}z)
  +b^{13}_-(q^{-k-2}z)-b^{23}_-(q^{-k-1}z)+(b+c)^{12}(q^{-k-1}z)}\Bigr):
  \label{F1}\\
  &&
  +f_{13}:e^{a^1_+(q^{\frac{k+1}{2}}z)
  -b^{23}_+(q^{k+1}z)-b^{13}(q^{k+1}z)+b^{23}(q^{k+2}z)}:, \n
  F^2(z)
  &\!\!=\!\!&
  \frac{1}{(q-q^{-1})z}:\Bigl(
  f_{21}e^{a^2_+(q^{\frac{k+1}{2}}z)-b^{23}(q^{k+1}z)}
  -f_{22}e^{a^2_-(q^{-\frac{k+1}{2}}z)-b^{23}(q^{-k-1}z)} \n
  &&\qquad
  -f_{23}e^{a^2_-(q^{-\frac{k+1}{2}}z)-b^{12}_-(q^{-k-1}z)
  -b^{13}_-(q^{-k-1}z)-(b+c)^{12}(q^{-k}z)-b^{13}(q^{-k}z)}
  \label{F2} \\
  &&\qquad
  +f_{24}e^{a^2_-(q^{-\frac{k+1}{2}}z)-b^{12}_+(q^{-k-1}z)
  -b^{13}_-(q^{-k-1}z)-(b+c)^{12}(q^{-k-2}z)-b^{13}(q^{-k}z)}
  \Bigr):,\nonumber
\ea
where $e_{11},e_{12},e_{21},e_{22}$,
$f_{11},f_{12},f_{13},f_{21},f_{22},f_{23},f_{24}$ are some constants.
{}From \eq{H1}\eq{H2}, we have
\ba
  H^1_n
  &\!\!=\!\!&
  a^1_nq^{-\frac12|n|}
  +b^{12}_nq^{-(\frac{k}{2}+1)|n|}(q^{|n|}+q^{-|n|})
  +b^{13}_nq^{-(\frac{k}{2}+2)|n|}
  -b^{23}_nq^{-(\frac{k}{2}+1)|n|},\\
  H^2_n
  &\!\!=\!\!&
  a^2_nq^{-\frac12|n|}
  -b^{12}_nq^{-(\frac{k}{2}+1)|n|}
  -b^{13}_nq^{-(\frac{k}{2}+1)|n|},\\
  K_1
  &\!\!=\!\!&
  q^{a^1_0+2b^{12}_0+b^{13}_0-b^{23}_0},\quad
  K_2=q^{a^2_0-b^{12}_0-b^{13}_0},\\  
  \psi^1_{\pm}(q^{\pm\frac{k}{2}}z)
  &\!\!=\!\!&
  e^{a^1_{\pm}(q^{\pm\frac{k+1}{2}}z)
  +b_{\pm}^{12}(q^{\pm k}z)+b_{\pm}^{12}(q^{\pm(k+2)}z)
  +b^{13}_{\pm}(q^{\pm(k+2)}z)-b^{23}_{\pm}(q^{\pm(k+1)}z)},\\
  \psi^2_{\pm}(q^{\pm\frac{k}{2}}z)
  &\!\!=\!\!&
  e^{a^2_{\pm}(q^{\pm\frac{k+1}{2}}z)
  -b^{12}_{\pm}(q^{\pm(k+1)}z)-b^{13}_{\pm}(q^{\pm(k+1)}z)}.
\ea
\begin{prop}
If $f_{11},f_{12},f_{13},f_{21},f_{22},f_{23},f_{24}$ satisfy
$$
  (f_{11},f_{12},f_{13},f_{21},f_{22},f_{23},f_{24})
  =
  (\frac{1}{e_{11}},\frac{1}{e_{12}},\frac{q^{k+1}e_{21}}{e_{11}e_{22}},
  \frac{q}{e_{21}},\frac{qe_{12}}{e_{11}e_{21}},\frac{1}{e_{22}},
  \frac{e_{12}}{e_{11}e_{22}}),
$$
then eqs. \eq{H1}-\eq{F2} realize $U_q(\widehat{sl}(2|1))$ with $\gamma=q^k$.
\end{prop}
{\it Proof:}
Straightforward OPE calculation shows this proposition.
For useful formulas, see ref.\cite{rAOS} (some modifications are
needed for $b^{ij}$ with $\nu_i\nu_j=-1$.). \hfill \qed

\noindent
{\bf Remark~}
Choices $e_{11}=e_{12}=e_{21}=e_{22}=1$ may be important.
Two of them, e.g. $e_{11}$ and $e_{21}$, is just a normalization of $F^i$'s.
When $e_{11}=e_{22}$, screening currents exist.

At present it is not clear that what choice of $e_{11},e_{12},e_{21},e_{22}$ 
is most natural. Perhaps it will be clarified by finding general formula for 
$U_q(\widehat{sl}(M|N))$. This problem is now under investigation.
For $M+N\leq 3$ realizations have been already obtained.
The Cartan part $H^i(z)$ for general $M$ and $N$ has also been 
obtained\cite{rAOSU} and we hope that we will be able to report a full answer.

\vskip5mm
\noindent{\bf Acknowledgments}\\
We would like to thank J.~Harvey, H.~Kubo, E.~Martinec, 
J.~Uchiyama and H.~Yamane
for valuable discussion.
We also would like to thank members of Nankai Institute for 
their hospitality since some of the results in section 3 were obtained
during our stay in Tianjin 19-24 August 1996.
This work is supported in part by Grant-in-Aid for Scientific
Research from Ministry of Science and Culture. 


\end{document}